\documentclass[a4paper,11pt,twocolumn,notitlepage]{article}

\newif\ifarxiv
\arxivtrue  

\usepackage[usenames,dvipsnames,svgnames,table]{xcolor}
\usepackage{bm}
\usepackage{amsmath,amsfonts,amssymb,mathtools}
\usepackage{graphicx}
\usepackage{nicefrac}
\usepackage{siunitx}
\usepackage{romannum}
\usepackage{enumitem}
\usepackage{algorithm}
\usepackage[noend]{algpseudocode}

\usepackage[left=.7in, right=.7in, top=1in, bottom=1in]{geometry}

\usepackage[
backend=biber,
style=ieee,
doi=false,
]{biblatex}
\defbibheading{refs}[\bibname]{%
	\section*{#1}%
	\markboth{#1}{#1}}

\usepackage{pgfplots}
\pgfplotsset{compat=1.18}
\usepackage{tikz}
\usepgfplotslibrary{fillbetween}
\usetikzlibrary{decorations.pathmorphing,
    patterns,shapes,arrows.meta,positioning, backgrounds, intersections, calc}
\tikzstyle{Nodes}=[circle, draw=black, fill=myblue!8, line width=1pt, minimum size=10pt]
\tikzstyle{Squares}=[rectangle, draw=black, fill=myblue!8, line width=1pt, minimum size=1pt]
\tikzstyle{Bars}=[rectangle, draw=black, fill=black, line width=.5pt, minimum size=1pt, inner sep=0pt, outer sep=0pt]
\tikzstyle{arrow} = [thick,-{Latex[length=2mm,width=2mm]}]
\tikzstyle{axis} = [very thick,-{Latex[length=2mm,width=2mm]}]
\usepackage{tkz-euclide}
\usetikzlibrary{pgfplots.groupplots}
\tikzstyle{block}=[
draw=white,
thick,
minimum width=1cm,
minimum height=1cm,
inner sep=4pt,
text width = 7em,
text centered,
fill = white,
]
\usetikzlibrary{external}
\usepackage{hyperref}
\hypersetup{
	colorlinks=true,
	citecolor=red,
    linkcolor=blue,
    filecolor=magenta,      
    urlcolor=blue,
	pdfauthor={Ruairi Moran},
	pdftitle={Paper},
	pdfkeywords={Paper, LaTeX},
	pdfsubject = {Paper},
}
\ifarxiv
    \newcommand{\pathtofigs}{./extern}
\else
    \newcommand{\pathtofigs}{../figures}
\fi

\newcommand{\mc}[1]{\mathcal{#1}}
\newcommand{\K}{\mc{K}}
\renewcommand{\S}{\mc{S}}

\newcommand\mymathbb[1]
{{\rm I\mathchoice{\hspace*{-2pt}}{\hspace*{-2pt}}
	{\hspace*{-1.75pt}}{\hspace*{-1.7pt}}#1}}  
\newcommand{\R}{\mymathbb R}				
\newcommand{\N}{\mymathbb N}
\newcommand{\E}{\mymathbb E}

\renewcommand{\P}{\mathbb{P}}
\newcommand{\smallplus}{{\scriptscriptstyle +}}
\newcommand{\Spp}{\mathbb{S}_{\smallplus\smallplus}}
\newcommand{\Sp}{\mathbb{S}_{\smallplus}}

\newcommand{\Rb}{\overline{\R}}

\newcommand{\xb}{\bar{x}}
\newcommand{\ub}{\bar{u}}

\newcommand{\zb}{\bar{z}}

\newcommand{\dom}{\operatorname{\mathbf{dom}}}
\newcommand{\diag}[1]{\operatorname{diag}\left(#1\right)}
\newcommand{\blkdiag}[1]{\operatorname{blkdiag}\left(#1\right)}
\newcommand{\range}{\operatorname{range}}
\newcommand{\rank}{\operatorname{\mathbf{rank}}}

\newcommand{\epi}{\operatorname{\mathbf{epi}}}

\newcommand{\prox}{\operatorname{\mathbf{prox}}}
\newcommand{\proj}{\operatorname{\mathbf{proj}}}

\newcommand{\AVaR}{\operatorname{AV@R}}
\newcommand{\VaR}{\operatorname{V@R}}

\newcommand{\SOC}{\operatorname{SOC}}
\newcommand{\indi}{\delta}

\DeclareMathOperator*{\argmin}{\mathbf{arg\,min}}

\DeclareMathOperator*{\minimise}{\mathbf{minimize}}

\newcommand{\tr}{\intercal}
\newcommand{\Tr}{^\tr}
\newcommand{\iTr}{^{i\tr}}

\newcommand{\imTr}{^{i_-\tr}}
\newcommand{\jTr}{^{j\tr}}

\let\oldmax\max
\renewcommand{\max}{\bm{\oldmax}}
\let\oldmin\min
\renewcommand{\min}{\bm{\oldmin}}
\let\oldinf\inf
\renewcommand{\inf}{\bm{\oldinf}}

\newcommand{\ri}{\mathbf{ri}}
\newcommand{\fix}{\mathbf{fix}\text{ }}
\newcommand{\id}{\mathrm{id}}

\newcommand{\smallmat}[1]{\left[ \begin{smallmatrix}#1 \end{smallmatrix} \right]}

\definecolor{myred}{rgb}{0.8,0.0,0.0}
\definecolor{mygreen}{rgb}{0.0,0.6,0.0}
\definecolor{myblue}{rgb}{0.0,0.0,0.8}

\newtheorem{theorem}{Theorem}

\newtheorem{definition}[theorem]{Definition}
\newtheorem{proposition}[theorem]{Proposition}

\newcommand{\nodes}{\operatorname{\mathbf{nodes}}}
\newcommand{\ch}{\operatorname{\mathbf{ch}}}
\newcommand{\anc}{\operatorname{\mathbf{anc}}}

\newcommand{\SPOCK}{\texttt{SPOCK}}
\newcommand{\MOSEK}{\texttt{MOSEK}}
\newcommand{\GUROBI}{\texttt{GUROBI}}

\newcommand{\IPOPT}{\texttt{IPOPT}}

\newcommand{\footremember}[2]{%
    \footnote{#2}
    \newcounter{#1}
    \setcounter{#1}{\value{footnote}}%
}
\newcommand{\footrecall}[1]{%
    \footnotemark[\value{#1}]%
} 

\title{GPU-Accelerated SPOCK for Scenario-Based Risk-Averse Optimal Control Problems\footnote{
This paper has received funding by the project
``BotDozer: GPU-accelerated model predictive control for autonomous heavy equipment,''
which is part of the Doctoral Training Grant No S3809ASA, funded by EPSRC.}}

\author{
Ruairi Moran\footremember{qub}{Queen's University Belfast, School of EEECS, i-AMS Centre, Ashby Building, BT9 5AH, Belfast, UK},
Pantelis Sopasakis\footrecall{qub}
}

\date{\today}

\begin{document}
\pagenumbering{arabic}

\maketitle

\begin{abstract}
    This paper presents a GPU-accelerated implementation of the \SPOCK{} algorithm, a proximal method designed for solving scenario-based risk-averse optimal control problems. The proposed implementation leverages the massive parallelization of the \SPOCK{} algorithm, and benchmarking against state-of-the-art interior-point solvers demonstrates GPU-accelerated \SPOCK's competitive execution time and memory footprint for large-scale problems. 
    We further investigate the effect of the scenario tree structure on parallelizability, and so on solve time.
\end{abstract}

\section{Introduction}
Risk-averse optimal control problems (RAOCPs) have a wide variety of applications, e.g., 
in MPC \cite{sopasakis_rampc} for microgrids \cite{hans_risk-measures}, 
in reinforcement learning \cite{chow_risk-measures2017, chow_risk-measures2015}, 
for energy management systems \cite{maree_risk-measures},
air-ground rendezvous \cite{barsi_risk-measures},
and collision avoidance \cite{dixit_risk-measures, schuurmans_risk-measures}.

In this work, the term \textit{risk} refers to the uncertainty of the detrimental effect of future events,
and the term \textit{risk-averse} is used to describe a preference for outcomes with reduced uncertainty of this detrimental effect.

The two main paradigms for approaching risk represent perspectives of the future. The \textit{stochastic/risk-neutral} paradigm represents the naive perspective that the future will unfold as expected. The \textit{robust/worst-case} paradigm represents the conservative perspective that the future will unfold in the worst possible way.
In terms of formulating optimal control problems (OCPs) for uncertain systems, the risk-neutral approach assumes that we have perfect information about the distributions of the involved random disturbances~\cite{cinquemani_stochastic_estimation}. The robust approach ignores any statistical information that is usually available~\cite[Ch. 3]{rawlings_mpc_book}.

The modern risk-averse paradigm forms a bridge between the risk-neutral and worst-case paradigms~\cite[Ch. 8]{dentcheva2024risk}.
The risk-averse formulation takes into account \textit{inexact} and \textit{data-driven} statistical information~\cite{wang_risk-problems-survey} using risk measures.
A \textit{risk measure} quantifies the magnitude of the right tail of a random cost, allowing the designer to strike a balance between the risk-neutral and the worst-case scenario.

While the ability to interpolate between the main paradigms is beneficial, the drawback is the computational burden; the cost functions of multistage RAOCPs involve the composition of several nonsmooth risk measures.

Historically, RAOCPs were (slowly) solved using stochastic dual dynamic programming (SDDP)~\cite{pereira_sddp-original, shapiro_analysis-sddp} approaches such as \cite{dacosta_sddp, shapiro_ra-and-rn-sddp}.
The modern approach~\cite{sopasakis2019ra-rc-mpc} reformulates RAOCPs to allow their solution through popular optimization software such as \MOSEK{}~\cite{andersen_mosek}, \GUROBI{}~\cite{gurobi_manual}, and \IPOPT{}~\cite{wachter2006implementation}.
The shortcomings of these interior-point solvers~\cite{nemirovski_interior-point} include scalability with the problem size
and significant memory requirements.

Despite the recent boom in GPU hardware and software, GPUs are mostly neglected in optimal control problems.
Successful examples in stochastic control are~\cite{chowdhury_gpu_stochastic_algo, sampathirao_stochastic_gpu}.

To address these issues, the \SPOCK{} algorithm~\cite[Algo. 3]{bodard2023spock} was recently proposed; this is a massively parallel proximal method for solving RAOCPs.
The \SPOCK{} algorithm employs the Chambolle-Pock (CP) method~\cite{chambolle_method} accelerated by the SuperMann framework~\cite{themelis_supermann} with Anderson's acceleration (AA)~\cite{anderson_acceleration_original}.
This allows the combination of the parallelisable nature of 
CP with the fast convergence properties of SuperMann.
Earlier, a serial implementation of \SPOCK{} in \texttt{Julia}~\cite{julia} was tested, and while the results were promising, the serial implementation ignored the rich structure of RAOCPs, leading to limited scalability.
Here we propose a GPU-accelerated implementation of \SPOCK{}, which exhibits high scalability, low solve times, and a low memory footprint.

The CP method was selected for its simple, gradient-free iterates.
While the CP method is similar to the celebrated alternating direction method of multipliers (ADMM)~\cite[Sec. 4.3]{chambolle_method}, ADMM requires an extra invocation of a computationally expensive linear operator. Indeed, our results found that the CP method outperforms ADMM for RAOCPs.

The main contributions of this paper are:
\begin{enumerate}
    \item a novel solver that does not suffer from the scalability issues of interior-point solvers,
    \item we demonstrate via extensive benchmarking that \SPOCK{} is amenable to parallelism on a GPU and compare 
    its performance against state-of-the-art solvers in terms 
    of computation time and memory usage,
    \item we propose a simple preconditioning of the problem data that leads to improved convergence properties, and
    \item we investigate the effect of the scenario tree structure on parallelizability, and so on solve time.
\end{enumerate}
The novel solver, GPU-accelerated \SPOCK{}, paves the way for real-time optimal control in high-dimensional and uncertain environments, offering a promising avenue for both academic research and practical applications.

\section{Notation}
Let $\N_{[k_1, k_2]}$ denote the integers in $[k_1, k_2]$,
let $1_n$ and $0_n$ refer to $n$-dimensional column vectors of ones and zeroes respectively, and
let $I_{n}$ be understood as the $n$-dimensional identity matrix.
We denote $n$-dimensional positive semidefinite and positive definite matrices by $\Sp^{n}$ and $\Spp^{n}$ respectively.
We denote the transpose of a matrix $A$ by $A\Tr$.
For $z\in\R^n$ let $[z]_+ = \max\{0,z\}$, where the max is taken element-wise.
The adjoint of a linear operator $L:\R^n\rightarrow\R^m$ is the operator $L^{\ast}:\R^m\rightarrow\R^n$ that satisfies $y\Tr Lx = x\Tr L^{\ast}y$ for all $x\in\R^n,y\in\R^m$.
The set of fixed points of an operator $T: \R^n \to \R^n$ is denoted $\fix T = \{v \in \R^n \mid v = T(v)\}$.
For two linear maps $F_1:\R^n\to\R^{m_1}$
and $F_2:\R^n\to\R^{m_2}$
we define $F_1 \times F_2:\R^n\to\R^{m_1}\times \R^{m_2}$
to be $(F_1\times F_2)(x) = (F_1(x), F_2(x))$.
For linear maps $F_i:\R^n\to\R^{m_i}$,
$i\in\N_{[1, K]}$,
we denote $\bigtimes_i F_i = F_1 \times \cdots \times F_K$.
For two linear maps $F_1:\R^{n_1}\to\R^{m_1}$ and $F_2:\R^{n_2}\to\R^{m_2}$, we define their direct sum as 
$F_1\oplus F_2\ni \R^{n_1+n_2}\ni(x_1, x_2)\mapsto (F_1(x_1), F_2(x_2))\in\R^{m_1+m_2}$.
The dual cone $\K^{\ast}$ of a closed convex cone $\K \subseteq \R^n$ is the set $\K^{\ast} = \{ y\in\R^n | y\Tr x \geq 0, \forall x\in\K \}$.
The relative interior of a convex set $\K\subseteq\R^n$ is denoted by $\ri\;\K$.
We denote the second order cone of dimension $d$ by $\SOC_d$ and a translated second order cone by $\SOC_d + a$ where $a\in\R^d$.
The indicator of a convex set $\K$ is defined as
$
    \indi_{\K}(x) = 0 \text{, if } x\in\K \text{, and }
    \indi_{\K}(x) = \infty \text{, otherwise. }
$
We define the set of extended reals,
$\Rb = \R \cup \{\infty\}$.
A function $f:\R^n\rightarrow\Rb$ is called \textit{lower semicontinuous} (lsc) if its sublevel sets, $\{ x \mid f(x) \leq \alpha \}$, are closed.
The domain of a function $f:\R^n\to\Rb$ is $\dom f = \{ x \mid f(x) < \infty \}$.
The subdifferential of a convex function $f:\R^n\to\Rb$ is $\partial f(x) = \{u\in\R^n \:|\: \forall y\in\R^n, (y-x)\Tr u + f(x) \leq f(y)\}$.
The conjugate of a function $f:\R^n\to\Rb$ is $f^{\ast}(y) = \sup_{x\in\R^n} \{ y\Tr x - f(x) \}$.
The proximal operator of a proper lsc convex function $f: \R^n \to \Rb$ with parameter $\gamma>0$ is
\begin{equation}
    \prox_{\gamma f}(x) = \argmin_v \left\{ f(v) + \tfrac{1}{2\gamma} \| v - x \|^2_2 \right\}.
\end{equation}
Where inequality signs are used with vectors (e.g., $x\leq y$, $x,y\in\R^n$), they are understood to hold element-wise.
Throughout, the notation superscript $i$ denotes scenario tree nodes, superscript $(k)$ denotes algorithm iterations, and subscript $t$ denotes time steps.

\section{Problem statement}
In this section we introduce multistage nested RAOCPs on scenario trees and define the problem statement.

\subsection{Dynamical system}
We consider the following discrete-time affine dynamical system
\begin{equation}\label{eq:disturbed-dt-affine-system}
    x_{t+1} = A(w_t)x_t + B(w_t)u_t + c(w_t),
\end{equation}
for $t\in\N_{[0,N-1]}$, with state variable $x_t\in\R^{n_x}$, control input $u_t\in\R^{n_u}$, and constant $c_t\in\R^{n_u}$. The random disturbance $w_t\in\{1,\ldots,n_w\}$ is finite valued; for example, it can be an iid process or a Markov chain.

\subsection{Scenario trees}\label{sec:scenario-trees-gpu}
\begin{figure}
    \centering
    \ifarxiv
        \includegraphics{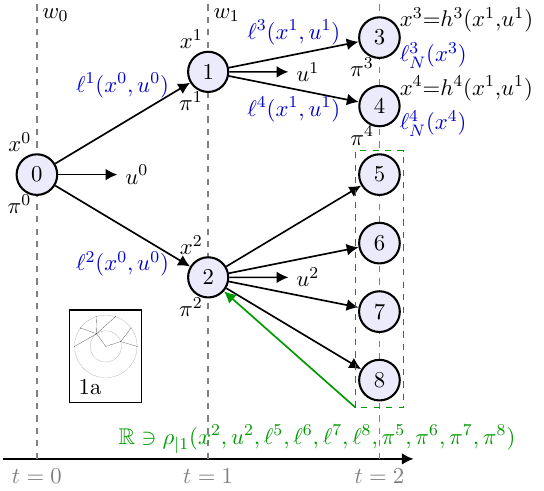}
    \else
        \input{\pathtofigs/trees.tex}
    \fi
    \caption{
        Discrete problem described by a scenario tree.
        Each node at stage $t+1$ is a possible realization of an event $w_t$ at stage $t$.
        The probability of each event is $\pi^i$.
        The dynamics $h^{i}$ and stage costs $\ell^{i}$ are a function of the ancestor state and input, and the event $w_t$ at the ancestor stage $t$.
        The terminal costs $\ell_{N}$ are a function of the ancestor state and the event $w_t$ at the ancestor stage $t$.
        The risk mapping $\rho$ is a many-to-one operator from children to their ancestor.
        Inset (1a) is a space efficient representation of the tree. The root node at stage $0$ is at the center and each consecutive stage is the next concentric circle outwards. This representation will be used later for depicting the branching structure of trees.
    }
    \label{fig:tree-example}
\end{figure}

A scenario tree is a representation of the dynamics of the system in Equation \eqref{eq:disturbed-dt-affine-system} over a finite number of stages given that the control actions are determined based on the system states in a causal fashion \cite{hoyland2001generating, dupavcova2000scenarios}, as depicted in Fig.~\ref{fig:tree-example}. Each time period is called a \textit{stage}, $t$, and the number of time periods is called the \textit{horizon}, $N$. For stopped processes, we denote the \textit{stopping} stage by $n_b$~\cite[Sec. 9.8]{gallager2013stochastic}. Each possible realization of Equation \eqref{eq:disturbed-dt-affine-system} at each stage is called a \textit{node}.
We enumerate the nodes of a tree with $i$, where $i=0$ is the \textit{root} node which corresponds to the initial state of the system. The nodes at subsequent stages for $t\in\N_{[0, N]}$ are denoted by $\nodes(t)$. Let $\nodes(t_1, t_2) = \cup_{t=t_1}^{t_2}\nodes(t)$, where $0\leq t_1\leq t_2\leq N$. The unique \textit{ancestor} of a node $i\in\nodes(1, N)$ is denoted by $\anc(i)$ and the set of \textit{children} of $i\in\nodes(t)$ for $t\in\N_{[0, N-1]}$ is $\ch(i) \subseteq \nodes(t+1)$. Each node is associated with a probability $\pi^i>0$ of occurring, where $\sum_{i\in\nodes(t)} \pi^i = 1$ for $t\in\N_{[0,N]}$.

\paragraph{Dynamics.}
The system dynamics $h^{i}$ on the scenario tree is described by
\begin{equation}\label{eq:finite-horizon-evolution-2-gpu}
    x^{i} = A^{i}x^{i_-} + B^{i}u^{i_-} + c^{i},
\end{equation}
where $i\in\nodes(0, N-1)$ and $i_-=\anc(i)$. The state $x^0$ at the root node is assumed to be known.

\paragraph{Constraints.}
Here we assume that the states and inputs must satisfy the random convex constraints
\begin{subequations}\label{eq:constraints-gpu}
    \begin{align}
        \Gamma_x^i x^i + \Gamma_u^i u^i \in{} & C^i,
        \\
        \Gamma_N^j x^j \in{}                  & C_N^j,
    \end{align}
\end{subequations}
for $i\in\nodes(0, N-1)$ and $j\in\nodes(N)$, where $C^i$ and $C_N^j$ are closed convex sets, $\Gamma_x^i$ and $\Gamma_u^i$ are the state-input constraint matrices and $\Gamma_N^j$ is the terminal state constraint matrix.

\paragraph{Quadratic-plus-linear costs.}
In this paper we consider quadratic-plus-linear cost functions.
On the scenario tree, for $i\in\nodes(0,N-1)$ and $i_{-}=\anc(i)$, every edge of the tree --- that connects the nodes $i$ and $i_{-}$ --- is associated with a stage cost function $\ell^{i}(x^{i_-}, u^{i_-}) = x\imTr Q^{i} x^{i_-} + u\imTr R^{i} u^{i_-} + q\iTr x^{i_-} + r\iTr u^{i_-}$.
For $j\in\nodes(N)$, every leaf node is associated with a terminal cost $\ell^{j}_N(x^j) = x\jTr Q^j_N x^j + q_N^{j\intercal}x^j$.

\subsection{Risk measures}\label{sec:risk-measures-gpu}
Given a discrete sample space $\Omega=\{1,\ldots, n\}$  and a probability vector $\pi\in\R^n$, a random cost $Z:\Omega\to\R$ can be identified by a vector $Z=(Z^1, \ldots, Z^n) \in \R^n$.
The expectation of $Z$ with respect to $\pi$ is $\E^{\pi}[Z] {}={} \pi\Tr Z$, and the maximum of $Z$ is defined as $\max Z = \max\{Z^i, i\in\N_{[1, n]}\}$.

A risk measure can be used to reflect the inexact knowledge of the probability measure.
A risk measure is an operator, $\rho:\R^n\to\R$, that maps a random cost $Z\in\R^n$ to a characteristic index $\rho(Z)$ that quantifies the magnitude of its right tail. Coherent risk measures~\cite[Sec. 6.3]{shapirolectures::2014} admit the following dual representation
\begin{equation}\label{eq:risk-dual-repr-gpu}
    \rho(Z) = \max_{\mu\in\mc{A}}\E^{\mu}[Z],
\end{equation}
where $\mc{A}$ is a nonemtpy closed convex 
subset of the probability simplex, 
$\Delta_n \coloneqq \{\mu\geq 0:\sum_i \mu_i = 1\}$,
known as the \textit{ambiguity set} of $\rho$. 
An ambiguity set allows us to quantify the 
uncertainty in $\pi$.
This means that coherent risk measures can be seen as 
a worst-case expectations of $Z$ 
with respect to a probability vector $\mu$ 
that is taken from $\mc{A}$. 
The expectation and the maximum can be thought of 
as the two extreme cases with $\mc{A}_{\E}=\{\pi\}$,
and $\mc{A}_{\max} = \Delta_n$.

\begin{figure}
    \centering
    \ifarxiv
        \includegraphics{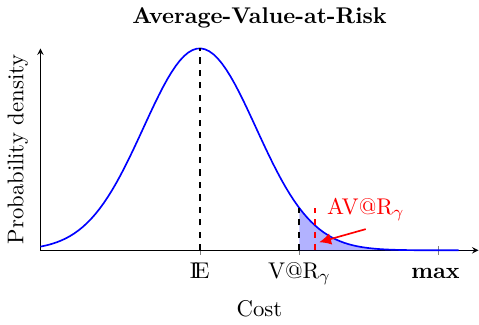}
    \else
        \input{\pathtofigs/avar.tex}
    \fi
    \caption{
    The value-at-risk at level $\gamma$ of a 
    continuously distributed random cost $Z$,
    denoted $\VaR_{\gamma}[Z]$,
    is the lowest cost $x$ such that $\mathrm{P}[Z\geq x]=\gamma$.
    The average-value at risk at level $\gamma$, $\AVaR_{\gamma}[Z]$,
    is the expectation of the tail of $Z$ 
    that lies above $\VaR_\gamma[Z]$~\cite{rockafellar2000optimization} (shaded area).%
    }
    \label{fig:avar}
\end{figure}

\paragraph{Average-value-at-risk.}
A popular coherent risk measure is the average value-at-risk with parameter $\gamma\in[0,1]$, denoted by $\AVaR_{\gamma}$.
It is defined as
\begin{equation}
    \AVaR_{\gamma}[Z] =
    \begin{cases}
        \min\limits_{t\in \R^n}
        \left\{
        t {+} \tfrac{1}{\gamma}\mathbb{E}^\pi[Z{-}t]_{+}
        \right\}, & \gamma {\neq} 0,
        \\
        \max[Z],  & \gamma {=} 0.
    \end{cases}
\end{equation}
The ambiguity set of $\AVaR_{\gamma}$ is
\begin{equation}
    \mathcal{A}_{\gamma}^{\rm avar}(\pi) = \left\{\mu {\in} \R^n \left| \sum_{i=1}^n \mu^i = 1, 0 \leq \mu^i \leq \tfrac{\pi^i}{\gamma} \right.\right\}.
\end{equation}
Figure~\ref{fig:avar} illustrates how $\AVaR_{\gamma}$ is an interpolation between the risk-neutral ($\AVaR_{1}=\E$) and worst-case ($\AVaR_{0}=\max$) approaches.

\paragraph{Conic representation.}
Coherent risk measures admit the following conic representation~\cite{sopasakis2019ra-rc-mpc}
\begin{equation}\label{eq:primal_risk-gpu}
    \rho[Z] = \max_{\mu \in \R^n, \nu \in \R^{n_{\nu}}}\
    \left\{ \mu\Tr Z \:|\: b - E\mu - F\nu \in \K \right\},
\end{equation}
where $\K$ is a closed convex cone.
Provided strong duality holds (i.e., if there exist $\mu$ and $\nu$ so $b - E\mu - F\nu \in \ri \K$), the dual form of the risk measure in Equation \eqref{eq:primal_risk-gpu} is
\begin{equation}\label{eq:dual_risk-gpu}
    \rho[Z] = \min_{y}\ \left\{ y\Tr b {}\mid{} E\Tr y = Z, F\Tr y = 0, y \in \K^* \right\}.
\end{equation}
For $\AVaR_{a}$, $n_{\nu}=0$,
$b =
    \begin{bmatrix}
        \pi\Tr & 0_n\Tr & 1
    \end{bmatrix}\Tr,$
$E =
    \begin{bmatrix}
        a I_n & -I_n & 1_n
    \end{bmatrix}\Tr,$
and
$\K = \R^{2n}_{+} \times \{0\}.$

\subsection{Conditional risk mappings}
Let $Z^i = \ell^{i}(x^{i_-}, u^{i_-})$ and $Z^j = \ell^{j}(x^j)$ for $i\in\nodes(1, N), j\in\nodes(N)$, and $i_-=\anc(i)$.
We define a random variable $Z_t = (Z^i)_{i\in\nodes(t+1)}$ on the probability space $\nodes(t+1)$, for $t\in\N_{[0,N-1]}$.
We define a random variable $Z_N = (Z^j)_{j\in\nodes(N)}$ on the probability space $\nodes(N)$.

We define $Z^{[i]} = (Z^{i_+})_{i_+\in\ch(i)}, i\in\nodes(0,N-1)$. This partitions the variable $Z_t = (Z^{[i]})_{i\in\nodes(t)}$ into groups of nodes which share a common ancestor.

Let $\rho: \R^{|\ch(i)|} \to \R$ be risk measures on the probability space $\ch(i)$. For every stage $t \in \N_{[0,N-1]}$ we may define a conditional risk mapping at stage $t$, $\rho_{|t}: \R^{|\nodes(t+1)|} \to \R^{|\nodes(t)|}$, as follows \cite[Sec. 6.8.2]{shapirolectures::2014}
\begin{equation}
    \rho_{|t}[Z_t] {}={} \left(\rho^i[Z^{[i]}]\right)_{i\in\nodes(t)}.
\end{equation}

\subsection{Original problem}
An RAOCP with horizon $N$ is defined via the following multistage nested formulation \cite[Sec. 6.8.1]{shapirolectures::2014}
\begin{multline}\label{eq:original-problem-gpu}
    V^{\star} {}={} \inf_{u_0} \rho_{|0} \bigg[ Z_0
        + \inf_{u_1} \rho_{|1} \Big[ Z_1
            + \cdots
            \\
            + \inf_{u_{N-1}} \rho_{|N-1} \big[ Z_{N-1}
                + Z_{N} \big] \cdots \Big]  \bigg],
\end{multline}
subject to \eqref{eq:finite-horizon-evolution-2-gpu} and \eqref{eq:constraints-gpu}.
Note that the infima in \eqref{eq:original-problem-gpu} are taken element-wise and with respect to the control laws $u_t = (u^i)_{i\in\nodes(t)}$ for $t\in\N_{[0, N-1]}$.
This problem is decomposed using \textit{risk-infimum interchangeability} and \textit{epigraphical relaxation} in \cite{sopasakis2019ra-rc-mpc}, allowing us to cast the original problem as the following minimization problem with conic constraints~\cite{bodard2023spock}
\begin{subequations}\label{eq:min-s0-v1-gpu}
    \begin{equation}
        \hspace*{-4.6cm}
        \minimise_{\substack{
                (x^{i})_{i}, (x^{j})_{j}, (u^{i})_{i}, (y^{i})_{i},
                \\
                (\tau^{i_+})_{i_+}, (s^{i})_{i}, (s^{j})_{j}}}\
        s^{0}
    \end{equation}
    \begin{align}
        \hspace*{-.2cm}
        \textbf{subject to}\
        & x^0 {}={} x,
        \\
        & x^{i} {}={} A^{i}x^{i_-} + B^{i}u^{i_-} + c^{i},
        \\
         & y^i \succcurlyeq_{(\K^i)^*} 0,
        y\iTr b^i {}\leq{} s^i,
        \\
         & E\iTr y^i {}={} \tau^{[i]} + s^{[i]},
        F\iTr y^i {}={} 0,
        \\
         & \ell^{i}(x^{i_-},u^{i_-}) \leq \tau^{i},
        \ell^{j}_N(x^{j}) \leq s^{j},
        \label{eq:cost-constraints-v1-gpu}
        \\
         & \Gamma_x^i x^i + \Gamma_u^i u^i \in C^i,
        \Gamma_N^j x^j \in C_N^j,
        \label{eq:constraint-2}
    \end{align}%
\end{subequations}%
for $t\in\N_{[0,N-1]}, i\in\nodes(t), i_+\in\ch(i), i_-=\anc(i)$, and $j\in\nodes(N)$, where $(\cdot)^{[i]} = ((\cdot)^{i_+})_{i_+\in\ch(i)}$, $\tau$ and $s$ are slack variables, and $x$ is the initial state.

\section{Reformulation}

\subsection{The Chambolle-Pock method}\label{sec:cp-method-gpu}
We can use the CP method to solve optimization problems of the following form \cite[p. 32]{ryu2016primer}
\begin{equation}\label{eq:main-cp-gpu}
    \mathbb{P}:\minimise_{z\in\R^{n_z}} f(z) + g(L z),
\end{equation}
where $L: \R^{n_z} \to \R^{n_{\eta}}$ is a linear operator, and $f$ and $g$ are proper closed convex functions on $\R^{n_z}$ and $\R^{n_{\eta}}$, respectively. We assume strong duality, that is, $\ri\,\textbf{dom} \: f \cap \ri\,\textbf{dom} \: g(L) \neq \emptyset$, throughout this work.

The CP method \cite{chambolle_method} recursively applies the firmly nonexpansive (FNE) operator $T$, where
$
    (
    z^{(k+1)}, \eta^{(k+1)}
    )
    =
    T(
    z^{(k)}, \eta^{(k)}
    ),
$
that is,
\begin{equation}\label{eq:chambolle-pock-operator-gpu}
    \begin{bmatrix}
        z^{(k+1)}
        \\
        \eta^{(k+1)}
    \end{bmatrix}
    =
    \underbracket[0.5pt]{
        \begin{bmatrix}
            \prox_{\alpha f} ( z^{(k)} {-} \alpha L^* \eta^{(k)} )
            \\
            \prox_{\alpha g^*} ( \eta^{(k)} {+} \alpha L (2 z^{(k+1)} {-} z^{(k)}) )
        \end{bmatrix}
    }_{T(z^{(k)}, \eta^{(k)})},
\end{equation}
which is a generalized proximal point method with preconditioning operator
\begin{equation}
    M(z, \eta) {}={}
    \begin{bmatrix}
        I         & -\alpha L^*
        \\
        -\alpha L & I
    \end{bmatrix}
    \begin{bmatrix}
        z
        \\
        \eta
    \end{bmatrix}.
\end{equation}
The CP method will converge to a solution, provided one exists and $0 < \alpha \|L\| < 1$, where the operator norm is $\|L\| = \max \left\{ \|L z\| : z\in\R^{n_{z}}, \|z\| \leq 1 \right\}$.

\subsection{Epigraphical relaxation of stage costs}
\label{sec:cost-constraint-reform}
In this section we show that the epigraphical relaxation of the stage costs in \eqref{eq:cost-constraints-v1-gpu} can be reformulated as second order conic (SOC) constraints.
The following proposition covers the common case where
$\ell$ is a convex quadratic-plus-linear function.
\begin{proposition}\label{prop:quad-plus-linear-to-soc}
    Let $Q\in\mathbb{S}^n_{+}$ with $\rank Q=p \leq n$,
    $q\in\R^n$ and define
    the function $\ell(z) = z^\intercal Q z + q^\intercal z$.
    Let $S\in\R^{n\times p}$ be such that
    $(\ker Q)^{\perp} = \range S$, so every $z\in\R^n$
    can be written uniquely as $z = z' + Se$, with
    $n\in\ker Q$ and $e\in\R^p$.
    Define $\tilde{Q} = S^{\intercal}QS$.
    Then, $(z, \tau)\in \epi \ell$ if and only if
    \begin{equation}
        \mathcal{G}_{S, Q, q}(e, z', \tau)
        \in
        \SOC_{p+2}
        +
        a_{S, Q, q},
        \label{eq:soc-condition-gpu}
    \end{equation}
    where $\mathcal{G}_{S, Q, q}$ is the linear map
    \begin{subequations}
        \begin{align}
            \mathcal{G}_{S, Q, q}(e, z', \tau) & {}={}
            \begin{bmatrix}
                \tilde{Q}^{1/2}e
                \\
                \tfrac{1}{2}\tau - \tfrac{1}{2}q^{\intercal}z'
                \\
                \tfrac{1}{2}\tau - \tfrac{1}{2}q^{\intercal}z'
            \end{bmatrix},
            \intertext{and}
            a_{Q, q}                     & {}={}
            \begin{bmatrix}
                - \tfrac{1}{2}\tilde{Q}^{-1/2}S^{\intercal}q
                \\
                - \tfrac{1}{8}\|S^\intercal q\|_{\tilde{Q}^{-1}}^2 + \tfrac{1}{2}
                \\
                - \tfrac{1}{8}\|S^\intercal q\|_{\tilde{Q}^{-1}}^2 - \tfrac{1}{2}
            \end{bmatrix}.
        \end{align}
    \end{subequations}
\end{proposition}
If $Q\in\mathbb{S}_{++}$, then $S = I$ and $\tilde{Q} = Q$;
in this case we will be using the notation $\mathcal{G}_{Q, q} = \mathcal{G}_{I, \tilde{Q}, q}$ and $a_{Q, q} = a_{I, \tilde{Q}, q}$.

\textit{Proof.}
We have $\ell(z) = z^\intercal Q z + q^\intercal z = r^\intercal \tilde{Q}r + q^\intercal S e + q^\intercal z'$.
By completing the squares we have
$\ell = \|e+\tfrac{1}{2}\tilde{Q}^{-1}S^{\intercal}q\|_{\tilde{Q}}^2 - \tfrac{1}{4}\|S^\intercal q\|_{\tilde{Q}^{-1}}^2$, therefore,
$(z, \tau)\in\epi \ell$ if and only if
\begin{align}
     & \left\|e + \tfrac{1}{2}\tilde{Q}^{-1}S^{\intercal}q\right\|_{\tilde{Q}}^2
    {}\leq{}
    \tau + \tfrac{1}{4}\|S^\intercal q\|_{\tilde{Q}^{-1}}^2
    \notag
    \\
    \Leftrightarrow
     &
    \left\|\tilde{Q}^{1/2}e + \tfrac{1}{2}\tilde{Q}^{-1/2}S^{\intercal}q\right\|_2^2
    {}\leq{}
    \tau + \tfrac{1}{4}\|S^\intercal q\|_{\tilde{Q}^{-1}}^2,
\end{align}
and using the identity
$\tau = \left(\frac{\tau+1}{2}\right)^2 - \left(\frac{\tau-1}{2}\right)^2$
we can show that this is equivalent to
the condition of Equation \eqref{eq:soc-condition-gpu}.~$\Box$

In some applications it is expedient to include an
$\ell_1$-type of penalty in the stage cost to enforce sparseness.
\begin{proposition}
    Let $Q\in\mathbb{S}^n_{++}$, $q\in\R^n$, $\lambda_0>0$, and define
    the function $\ell(z) = z^\intercal Q z + q^\intercal z + \lambda_0 \|z\|_1$.
    Then $(z, \tau) \in \epi \ell$ if and only if there exists scalars
    $\theta\in\R$ and $\lambda_1, \ldots, \lambda_n \geq 0$ such that
    \begin{equation}
        \mathcal{G}_{Q, q}(e, z', \tau-\lambda_0\theta)
        \in
        \SOC_{p+2}
        +
        a_{Q, q},
        \label{eq:soc-condition-with-norm-1-gpu}
    \end{equation}
    and $\sum_{i=1}^{n}\lambda_i \leq \theta$ and $-\lambda_i \leq z_i \leq \lambda_i$
    for $i=1,\ldots, n$.
\end{proposition}
\textit{Proof.} We have that $z^\intercal Q z + q^\intercal z + \lambda_0 \|z\|_1 \leq \tau$ is
equivalent to the existence of a $\theta\in\R$ such that
\begin{subequations}
    \begin{align}
        z^\intercal Q z + q^\intercal z + \lambda_0 \theta \leq \tau,
        \\
        \|z\|_1 \leq \theta,
    \end{align}
\end{subequations}
where, by Proposition \ref{prop:quad-plus-linear-to-soc},
the first condition yields Equation \eqref{eq:soc-condition-with-norm-1-gpu}.
The second condition is equivalent to
$\sum_{i=1}^{n}|z_i| \leq \theta$
and using the relaxation $|z_i|\leq \theta$
completes the proof.~$\Box$

Following the same approach we can accommodate 
stage cost functions that involve soft constraints 
with quadratic penalty functions, e.g., 
functions of the form 
$\ell(z) = z^\intercal Q z + q^\intercal z
+ \|\max(0, z-z_{\rm max})\|^2$, 
with $Q\in\mathbb{S}_{+}^{n}$,
where the maximum 
is meant in the element-wise sense.
The constraint $\ell(z)\leq \tau$ can be 
written as
\begin{equation}
    \smallmat{z \\ \theta}^\intercal 
    \smallmat{Q \\ & 1}
    \smallmat{z \\ \theta}
    +
    \smallmat{q \\ 0}^\intercal
    \smallmat{z \\ \theta}
    \leq 
    \tau,
    \label{eq:soft-constraint-relaxation}
\end{equation}
with $\theta \geq 0$ and $x - x_{\rm max} \leq \theta$,
and Equation \eqref{eq:soft-constraint-relaxation}
can be handled using 
Proposition \ref{prop:quad-plus-linear-to-soc}.

\subsection{Sets}
\label{sec:sets}
Define $z = (s^0, z_1, z_2)$, where
$
    z_1 = ((x^{i'})_{i'}, (u^i)_{i}),
$ and $
    z_2 = (y^i,\tau^{[i]},s^{[i]}),
$
for $i'\in \nodes(0, N)$ and $i \in \nodes(0, N-1)$.
We define the sets
\begin{align}
    \S_1
    =          &
    \left\{
    z_1
    \left|
    \begin{array}{l}
        x^0 - x = 0,
        \\
        x^{i} - A^{i} x^{i_-} - B^{i} u^{i_-} - c^{i} = 0,
        \\
        \forall i\in\nodes(1,N), i_-=\anc(i).
    \end{array}
    \right.
    \right\},
    \\
    \S_2
    =          & \prod_{i\in\nodes{(0,N-1)}}
    \underbracket[0.5pt]{
        \ker
        \begin{bmatrix}
            E\iTr & -I & -I
            \\
            F\iTr & 0  & 0
        \end{bmatrix}
    }_{\S_2^i},
    \\ \label{eq:S3-def-gpu}
    \S_3
    =          &
    \prod_{i\in\nodes{(0,N-1)}}
    \left((\K^i)^* \times \R_+ \times C^i\right)
    \notag
    \\
    {}\times{} &
    \prod_{i\in\nodes{(1,N)}}
    \left( \SOC_{n_x + n_u + 2} + a^i\right)
    \notag
    \\
    {}\times{} &
    \hspace*{-1.2em}\prod_{j\in\nodes{(N)}}
    \left(C^j
    \times
    \hspace*{-0.2em}\left(
        \SOC_{n_x + 2} + a_N^j \right)
    \right),
\end{align}
where $a^i = a_{\blkdiag{Q^i, R^i}, [q^{i\intercal}~ r^{i\intercal}]^\intercal}$ and $a_N^j = a_{Q_N^j, q_N^j}$.

\subsection{Linear operator $L$}\label{sec:linear-operator}
We define the linear operator $L$ that maps $z$ to
\begin{align}
    \eta =
    (
     & (
    y^i,
    s^i - b\iTr y^i,
    \Gamma_x^{i} x^i + \Gamma_u^{i} u^i
    )_{i\in\nodes(0,N-1)},
    \notag
    \\
     & (
    ( Q^i )^{\nicefrac{1}{2}} x^{\anc(i)},
    ( R^i )^{\nicefrac{1}{2}} u^{\anc(i)},
    \notag \\ & \qquad\qquad\qquad  
    \tfrac{1}{2}\tau^{i},
    \tfrac{1}{2}\tau^{i}
    )_{i\in\nodes(1, N)},
    \notag
    \\
     & (
    \Gamma_N^j x^j,
    (( Q_N^j )^{\nicefrac{1}{2}} x^j,
    \tfrac{1}{2}s^j,
    \tfrac{1}{2}s^j)
    )_{j\in\nodes(N)}
    ).
\end{align}

It can be seen that there is a permutation 
$\Pi: \R^{n_z}\to\R^{n_z}$ such that $L\circ \Pi$
is block-diagonal with 
\begin{equation}
    L \circ \Pi 
    {}={} 
    \hspace{-1.5em}
    \bigoplus_{i\in\nodes(0, N-1)}
    \hspace{-1em}
    L_i
    \oplus 
    \hspace{-0.75em}\bigoplus_{j\in\nodes(N)}\hspace{-1em} 
    {}L_N^j{}
    \oplus 
    D_0,
\end{equation} 
where
\begin{align}
    D_i
    {}={} &
    \bigtimes_{i_+\in\ch(i)}\blkdiag{Q^{i_+}, R^{i_+}, \tfrac{1}{2}1_{2}}
    \\
    D^j_N
    {}={} &
    \blkdiag{Q^j, \tfrac{1}{2}1_{2}}
    \\
    b
    {}={} &
    \hspace{-1em}\bigtimes_{i\in\nodes(0,N-1)}\hspace{-1em}b^i,
    \\
    D_0
    {}={} &
    \begin{bmatrix}
        1\Tr_{n_b} & 0\Tr_{|\nodes(0,N-1)|}
        \\
        b^*        & 1\Tr_{|\nodes(0,N-1)|}
    \end{bmatrix},
    \\
    L^i
    {}={} &
    D^i
    \times
    \begin{bmatrix}
        \Gamma_x^i & \Gamma_u^i & 0_{n_g \times n_c}
    \end{bmatrix},
    \\
    L^j_N
    {}={} &
    D^j_N
    \times
    \begin{bmatrix}
        \Gamma_N^j & 0_{n_g}
    \end{bmatrix}.
\end{align}
It follows that 
\begin{equation}
    \|L\| = \|L\circ \Pi\| = \max\left\{ (\|L_i\|)_i, (\|L_N^j\|)_j, \|D_0\| \right\}.
\end{equation}
Using the property $\|\bigtimes_{i=1}^{\upsilon} L_i\| \leq \sqrt{\upsilon}\max_{i=1,\ldots, \upsilon}\|L_i\|$, it follows that $\|L\|$ is controlled by 
\begin{multline}
    \max_{i, j}\{\sqrt{|\ch(i)|}\max(\|Q^{i_+}\|, \|R^{i_+}\|), 
    \\
    \|\smallmat{\Gamma_x^i & \Gamma_u^i}\|,
    \|Q^j\|,
    \|\Gamma_N^j\|\}
\end{multline}
In other words, $\|L\|$---and, importantly, the CP step size $\alpha$---scales with the largest number of children per node, not the prediction horizon $N$.

\subsection{Operator splitting}
\label{sec:operator-splitting-gpu}
We can now write Problem \eqref{eq:min-s0-v1-gpu} as
\begin{subequations}
    \begin{align}
        \minimise_{z}\  & s^0
        \\
        \textbf{subject to}\
                        & z_1 \in \S_1,
        z_2 \in \S_2,
        \eta \in \S_3,
    \end{align}
\end{subequations}
which is equivalent to \eqref{eq:main-cp-gpu}
\begin{equation}
    \minimise_{z}\ f(z) + g(Lz),
    \label{eq:cp-general-splitting-gpu}
\end{equation}
where
\begin{subequations}
    \begin{align}
         & f(z) {}={}
        s^0
        {}+{}
        \delta_{\S_1}(z_1)
        {}+{}
        \delta_{\S_2}(z_2),
        \label{eq:f-gpu}
        \\
         & g(\eta) {}={} \delta_{\S_3}(\eta).
        \label{eq:g-gpu}
    \end{align}
\end{subequations}

\section{Numerical algorithm}

\subsection{Proximal operators}
The proximal operators in \eqref{eq:chambolle-pock-operator-gpu} remain largely the same as in~\cite{bodard2023spock}, however, here we extend the projection on the dynamics to allow the affine set $\S_1$. We denote the primal vector before projection by $\zb$. Computing $\proj_{\S_1}(\zb_1)$ comprises the offline Algorithm~\ref{alg:s1projection:offline-gpu} from~\cite{bodard2023spock}, which we show here for completeness, and the online Algorithm~\ref{alg:s1projection:online-gpu}. Algorithm~\ref{alg:s1projection:offline-gpu} is run offline once, and Algorithm~\ref{alg:s1projection:online-gpu} is run online once per projection on $\S_1$.
\begin{algorithm}[htbp!]
    \caption{Projection on $\S_1$: Offline}\label{alg:s1projection:offline-gpu}
    \begin{algorithmic}[1]
        \Require Matrices $(A^i)_i$ and $(B^i)_i$, and horizon $N$

        \Ensure Matrices $(K^i)_i$, $(P^i)_i$, $(\widetilde{R}^i)_i$ and $(\bar{A}^i)_i$

        \ForAll{$i\in\nodes(N)$ \textbf{in parallel}}
        \State $P^i \gets I_{n_x}$
        \EndFor

        \For{$t=0,1,\ldots, N-1$}
        \ForAll{$i{\in}\nodes(N{-}(t{+}1)) \textbf{ in parallel}$}
        \State$
            \widetilde{P}^{i_+} \gets
            B^{i_+\intercal} P^{i_+}
        $, $i_+{\in}\ch(i)$
        \State$
            \widetilde{R}^i \gets $
            $
                I_{n_u}
                + \sum_{i_+\in\ch(i)} \widetilde{P}^{i_+} B^{i_+}
            $
            and determine its Cholesky decomposition
        \State$
            K^i \gets
            -(\widetilde{R}^{i})^{-1}
            \sum_{i_+\in\ch(i)} \widetilde{P}^{i_+} A^{i_+}
        $
        \State$
            \bar{A}^{i_+} \gets
            A^{i_+} + B^{i_+}K^i
        $, $i_+{\in}\ch(i)$
        \State$
            P^i \gets
            I_{n_x} {+} K^{i\intercal} K^i {+} \sum_{i_+\in\ch(i)} \bar{A}^{i_+\intercal} P^{i_+} \bar{A}^{i_+}
        $
        \EndFor

        \EndFor

    \end{algorithmic}
\end{algorithm}
\begin{algorithm}[htbp!]
    \caption{Projection on $\S_1$: Online}\label{alg:s1projection:online-gpu}
    \begin{algorithmic}[1]
        \Require Matrices computed by Algorithm \ref{alg:s1projection:offline-gpu},
        constants $(c^i)_i$, vector $\zb_1$ as defined in Sec.~\ref{sec:operator-splitting-gpu}, and initial state $x^0$

        \Ensure Projection on $\S_1$ at $\zb_1$

        \ForAll{$i\in\nodes(N)$ \textbf{in parallel}}
        \State $q^i_0 \gets -\xb^i$
        \EndFor

        \For{$t=0,1,\ldots, N-1$}
        \ForAll{$i{\in}\nodes(N{-}(t{+}1)) \textbf{ in parallel}$}
        \State$
            \bar{d}^i \gets
            \sum_{i_+\in\ch(i)} B^{i_+\intercal}
            \left(q^{i_+} + P^{i_+} c^{i_+} \right)
        $
        \State$
            d_{t+1}^i \gets
                (\widetilde{R}^{i})^{-1}
                \Big(
                \ub^i -
                \bar{d}^i
                \Big)
            $
            using Cholesky of $\widetilde{R}^{i}$
        \State$
            \bar{q}^i \gets
                \sum_{i_+\in\ch(i)} \bar{A}^{i_+\intercal} 
                ( 
                P^{i_+}
                ( B^{i_+} d^i
                + c^{i_+} )
                + q^{i_+} 
                )
            $
        \State$
            q_{t+1}^i \gets
            K^{i\intercal} \left( d^i - \ub^i \right)
            - \xb^i +
            \bar{q}^i
            $
        \EndFor
        \EndFor

        \State $x^0 \gets x$
        \For{$t=0,1,\ldots, N-1$}
        \ForAll{$i{\in}\nodes(t),i_+{\in}\ch(i)$ \textbf{in parallel}}
        \State $u^i \gets K^i x^i + d^i$
        \State $x^{i_+} \gets A^{i_+}x^i + B^{i_+}u^i + c^{i_+}$
        \EndFor
        \EndFor

        \State \textbf{return} $z_1 \gets \proj_{\S_1}(\zb_1)$.

    \end{algorithmic}
\end{algorithm}

\subsection{\SPOCK{} algorithm}
We state the parallelized \SPOCK{} algorithm from \cite[Algo. 3]{bodard2023spock}. We define 
$
v = 
\begin{bmatrix}
    z\Tr & \eta\Tr
\end{bmatrix}\Tr
$. 
Note that we do not explicitly form or store the matrices $L, L^{\ast}$ or $M$.
The \SPOCK{} algorithm follows the SuperMann~\cite{themelis_supermann} framework, a Newton-type algorithm that finds a fixed point of firmly nonexpansive operators. 
The framework involves two extragradient-type updates that can involve fast (e.g., quasi-Newtonian) directions to enjoy the same global convergence properties as the classical Krasnosel'ski\v{\i}-Mann scheme.

These directions are computed using Anderson's acceleration~\cite{anderson_acceleration_original} as described in Algorithm~\ref{alg:aa}, and is discussed further in Section~\ref{sec:aa-parallel} along with the choice of parameter $m$.
AA uses two matrices, $M_{r}, M_{d} \in\R^{n_v \times m}$, to keep a running store of the last $m$ residuals and their differences, respectively.
\begin{algorithm}
    \caption{Anderson's acceleration}\label{alg:aa}
    \begin{algorithmic}[1]
        \Require
        Memory length $m$, residual $r^{(k)}$, and matrices $M_{r}$ and $M_{d}$.

        \Ensure Fast direction $\psi^{(k)}$ and matrices $M_{r}$ and $M_{d}$

        \State Shift columns of $M_{r}$ and $M_{d}$ to the right by $1$

        \State $M_{r}(:,0) \gets r^{(k)}$
        
        \State $M_{d}(:,0) \gets M_{r}(:,0) - M_{r}(:,1)$

        \If{$k \leq m$}
        \State $\psi^{(k)} \gets -r^{(k)}$
        \Return $\psi^{(k)}$
        \EndIf

        \If{columns of $M_{r}$ and $M_{d}$ are greater than $m$}
        \State Delete last column of $M_{r}$ and $M_{d}$
        \EndIf

        \State $Q^{(k)},R^{(k)} \gets M_{d}$ by QR decomposition

        \State Solve $R^{(k)} \kappa = Q^{(k)\intercal} r^{(k)}$ for $\kappa$

        \State $\psi^{(k)} \gets -r^{(k)} - (M_{r} - M_{d})\kappa$
        
        \State \textbf{return} $\psi^{(k)}$
    \end{algorithmic}
\end{algorithm}

For the other \SPOCK{} parameters, we refer to~\cite[Sec. VI.D]{themelis_supermann}, that is, set 
$c_0 = 0.99$, 
$c_1 = 0.99$,
$c_2 = 0.99$,
$\beta = 0.5$, 
$\sigma = 0.1$, 
and $\lambda = 1$.
Indeed, these values have shown to be effective after some parameter testing.
\begin{algorithm}
    \caption{\SPOCK{} algorithm for RAOCPs}\label{alg:spock-gpu}
    \begin{algorithmic}[1]
        \Require
        problem data,
        $z^{(0)}$ and $\eta^{(0)}$,
        tolerances $\epsilon_{\text{abs}} > 0$ and $\epsilon_{\text{rel}} > 0$,
        $\alpha$ such that $0 < \alpha \|L\| < 1$,
        $\N \ni m > 0$,
        $c_0, c_1, c_2 \in [0,1)$,
        $\beta,\sigma \in (0,1)$, and
        $\lambda \in (0,2)$.

        \Ensure approximate solution of RAOCP

        \State Compute offline matrices via Algorithm \ref{alg:s1projection:offline-gpu}

        \State $r^{(0)} \gets v^{(0)} - T(v^{(0)})$,
        $\zeta_{(0)} {\gets} \|r^{(0)}\|_M$

        \State $\omega_{\text{safe}} \gets \zeta_{(0)}$, $k \gets 0$

        \If{termination criteria (Sec.~\ref{sec:termination-conditions}) are satisfied}\label{line:check-term-gpu}
        \State \textbf{return} $v^{(k)}$
        \EndIf

        \State $r^{(k)} \gets v^{(k)} - T(v^{(k)})$

        \State $\psi^{(k)} \gets$ Algorithm~\ref{alg:aa}\label{line:aa}

        \State $\omega_{(k)} {\gets} \|r^{(k)}\|_M$

        \If{$\omega_{(k)} \leq c_0 \zeta_{(k)}$}
        \State $v^{(k+1)} {\gets} v^{(k)} {+} \psi^{(k)}$,
        \hfill{\color{gray}(\textit{K0})}

        $\zeta_{(k+1)} {\gets} \omega_{(k)}$, \textbf{goto} line \ref{line:last-step-gpu}
        \EndIf

        \State $\zeta_{(k+1)} \gets \zeta_{(k)}$, $\tau \gets 1$

        \State $\tilde{v}^{(k)} \gets v^{(k)} + \tau \psi^{(k)}$, $\tilde{r}^{(k)} \gets \tilde{v}^{(k)} - T(\tilde{v}^{(k)})$ \label{line:try-new-tau-gpu}

        \State $\tilde{\omega}_{(k)} {\gets} \|\tilde{r}^{(k)}\|_M$

        \If{$\omega_{(k)} \leq \omega_{\text{safe}}$ \textbf{and} $\tilde{\omega}_{(k)} \leq c_1 \omega_{(k)}$}
        \State $v^{(k+1)} {\gets} \tilde{v}^{(k)}$,
        \hfill{\color{gray}(\textit{K1})}

        $\omega_{\text{safe}} {\gets} \tilde{\omega}_{(k)} {+} c_2^k$, \textbf{goto} line \ref{line:last-step-gpu}
        \EndIf

        \State $\rho_{(k)} \gets \tilde{\omega}_{(k)}^2 - 2\alpha (\tilde{r}^{(k)})\Tr M ( \tilde{v}^{(k)} - v^{(k)} )$

        \If{$\rho_{(k)} \geq \sigma \tilde{\omega}_{(k)} \omega_{(k)}$}
        \State $v^{(k+1)} {\gets} v^{(k)} {-} \dfrac{\lambda \rho_{(k)}}{\tilde{\omega}_{(k)}^2} \tilde{r}^{(k)}$ \hfill{\color{gray}(\textit{K2})}
        \EndIf
        \State \textbf{else} $\tau \gets \beta \tau$, \textbf{goto} line \ref{line:try-new-tau-gpu}

        \State $k \gets k + 1$, \textbf{goto} line \ref{line:check-term-gpu} \label{line:last-step-gpu}
    \end{algorithmic}
\end{algorithm}

\subsection{Termination conditions}\label{sec:termination-conditions}
We define the residual operator $R = \id - T$, where $\id$ is the identity operator.
The SuperMann algorithm finds a fixed point $v^{\star}\in\fix T$ of a firmly nonexpansive operator $T$ by finding a zero of $R$. Therefore, it makes sense to terminate \SPOCK{} when $\|R^{(k)}\|_{\infty} \leq \max(\epsilon_{\rm abs}, \epsilon_{\rm rel}\|R^{(0)}\|_{\infty})$ for some small $0 < \epsilon_{\rm abs}, \epsilon_{\rm rel} \in\R$. However, this does not guarantee $v^{\star}$ is an $\epsilon$-approximate optimal point to $\P$, as we now discuss.

The optimality condition for $\P$ is
\(
0 \in \partial f(z) + L^* \partial g(L z)
\)
assuming $\ri\,\textbf{dom} \: f {\cap} \ri\,\textbf{dom} \: g(L) {\neq} \emptyset$, or
\begin{subequations}%
    \begin{align}%
        0 & \in \partial f(z) + L^* \eta,
        \label{eq:primal-dual-point:a-gpu}
        \\
        0 & \in \partial g^*(\eta) - Lz,
        \label{eq:primal-dual-point:b-gpu}
    \end{align}
\end{subequations}
for some $z\in\R^n$ and $\eta \in \R^m$.

To define a notion of approximate optimality, we firstly
introduce a notion of $\epsilon$-approximate belonging to a set.
For a set $\mc{X}\subseteq\R^n$, the notation $x \in_{\epsilon} \mc{X}$
means that there is a $\xi\in\R^n$, with $\|\xi\|_\infty\leq \epsilon$ and $x+\xi\in \mc{X}$.
\begin{definition}[Approximate optimality]
    A point $(\tilde{z}, \tilde{\eta})$ is said to be an
    $\epsilon$-approximate primal-dual optimal point for
    $\mathbb{P}$ if
    $0 \in_\epsilon{} \partial f(\tilde{z}) + L^* \tilde{\eta}$
    (primal condition)
    and
    $0 \in_\epsilon{} \partial g^*(\tilde{\eta}) - L \tilde{z}$
    (dual condition).
\end{definition}

Next, we turn our attention to the CP operator to
determine appropriate primal and dual residuals.
The primal update in \eqref{eq:chambolle-pock-operator-gpu}
can be written as
\begin{align}
                      & \alpha \partial f ( z^{(k+1)} ) + z^{{k+1}} \ni z^{(k)} - \alpha L^* \eta^{(k)}
    \notag
    \\
    \Leftrightarrow\: &
    \underbracket[0.5pt]{\tfrac{\Delta z^k}{\alpha} - L^*(\Delta \eta^k)}_{\xi_1^{(k)}} \in L^*\eta^{(k+1)} + \partial f(z^{(k+1)}),
    \label{eq:xi1k}
\end{align}
with $\Delta z^{(k)} = z^{(k)} - z^{(k+1)}$ and
$\Delta \eta^{(k)} = \eta^{(k)} - \eta^{(k+1)}$.
We, therefore, define the residual $\xi_1^{(k)}$
as in Equation \eqref{eq:xi1k}.
Likewise, we define the residual $\xi_2^{(k)}$
\begin{align}
    & \alpha \partial g^* (\eta^{(k+1)}) {+} \eta^{(k+1)} {-} 2\alpha L z^{(k+1)} \ni \eta^k {-} \alpha L z^{(k)}
    \notag
    \\
    \Leftrightarrow\: &
    \underbracket[0.5pt]{\tfrac{\Delta\eta^k}{\alpha} - L(\Delta z^{(k)})}_{\xi_2^{(k)}}
    \in \partial g^* (\eta^{(k+1)}) - L z^{(k+1)}.
\end{align}
If $\|\xi_1^{(k+1)}\|_\infty \leq \epsilon$,
and $\|\xi_2^{(k+1)}\|_\infty \leq \epsilon$,
then $(z^{(k+1)}, \eta^{(k+1)})$ is an $\epsilon$-approximate optimal point.

Remark: suppose $(z^{(k+1)}, \eta^{(k+1)})$
is an $\epsilon$-approximate optimal point.
Then, $Lz^{(k+1)} \in_\epsilon \partial g^*(\eta^{(k+1)})$.
Let $(z^\star, \eta^\star)$ be a fixed point of the CP method.
Then, $Lz^\star \in \partial g^*(\eta^\star)$.
If $\partial g^*$ is metrically regular
with modulus $\tau$
at $(\eta^\star, L z^\star)$~\cite[Sect.~9.G]{rockafellarWets2009},
then $\eta^{(k+1)} \in_{\epsilon \tau} \partial g(Lz^{(k+1)})$.
As a result,
$0 {}\in_{\epsilon'}{} \partial f(z^{(k+1)}) + L^* \partial g(Lz^{(k+1)})$ with $\epsilon' = \epsilon (1 + \tau\|L^*\|_{\infty})$,
i.e., the pair $(z^{(k+1)}, \eta^{(k+1)})$
satisfies approximately the original
optimality conditions.

\subsection{Preconditioning}
Both CP and SPOCK are first-order methods, 
so they can benefit from preconditioning. 
It is possible to precondition the CP method
at the expense of a higher per-iteration computation 
cost~\cite{chambolle_preconditioning}.
Here, however, we use heuristic preconditioning matrices to scale the problem --- that is, increase the step size --- offline.
In Section~\ref{sec:linear-operator}, we see that the step size of the CP method is proportional 
to $\|L\|$, and so to the cost matrices $Q$ and $R$.
As a result, we scale the states and inputs by the diagonal matrices $S^{x}\in\Spp^{n_x}$ and $S^{u}\in\Spp^{n_u}$ respectively, where
\begin{subequations}
\begin{align}
    \hat{c}
    {}={}&
    \sqrt{\max(|\ch(i)|)}
    \\
    S^{x}_{kk} 
    {}={}& 
    \hat{c}
    \max\left\{
        1,
        ((Q^{i}_{kk})^{\nicefrac{1}{2}})_{i}
    \right\},
    \\
    S^{u}_{kk} 
    {}={}& 
    \hat{c}
    \max\left\{
        1,
        ((R^{i}_{kk})^{\nicefrac{1}{2}})_{i}
    \right\},
    \\
    S^{x}_{N,kk} 
    {}={}& 
    \max\left\{
        1,
        ((Q^{j}_{N,kk})^{\nicefrac{1}{2}})_{j}
    \right\},
\end{align}
\end{subequations}
for $i\in\nodes(1,N)$ and $j\in\nodes(N)$.
In other words, we introduce a change of variables 
$\tilde{x}^i = S^x x^i$, and  $\tilde{u}^i = S^u u^i$
for $i\in\nodes(1, N-1)$,
and $\tilde{x}^j = S^x_N x^j$, $j\in\nodes(1, N)$.
Alongside, the tolerance in the termination criteria is 
scaled accordingly.
The constraints become 
$\Gamma_x^i S_x^{-1}\tilde{x}^i + \Gamma_u^i S_u^{-1} \tilde{u}^i \in C^i$ and $\Gamma_N^j S_x \tilde{x}^j \in C_N^j$
for $i\in\nodes(0,N-1)$ and $j\in\nodes(N)$.
We define $a^i = \max\{1, \|\smallmat{\Gamma_x^i S_x^{-1} &  \Gamma_u^i S_u^{-1}}\|\}$ and scale the constraints as 
$\tfrac{1}{a^i}\Gamma_x^i S_x^{-1}\tilde{x}^i + \tfrac{1}{a^i}\Gamma_u^i S_u^{-1} \tilde{u}^i \in \tfrac{1}{a^i}C^i$.

This change of variables and rescaling of the constraints 
is shown to prevent $\|L\|$ 
from becoming too large and the step size from becoming too small. 

\section{Parallel computing}
There are two important factors, data type and block size, that must be considered when dealing with parallel computing. There is a trade-off between speed and accuracy in data types, that is, using the \texttt{float} type will decrease execution time at the cost of accuracy compared to using the \texttt{double} type. The block size, or threads per block, is a programming abstraction that determines how many threads are executed in parallel. While choosing the optimal block size is hardware specific, it is recommended to use multiples of the warp size~\cite{practices2020cuda} (32 for all current NVIDIA hardware and compute capabilities~\cite[Tab. 16]{programming2020cuda}) so that warps are fully utilized.

For \SPOCK{} in particular, we recommend the \texttt{float} type for solving problems to lax tolerances ($10^{-3}$), and the \texttt{double} type for strict tolerances ($10^{-5}$). Furthermore, we recommend the block size to be at least the number of states plus the number of inputs plus two, i.e., $n_x + n_u + 2$. This ensures optimal performance of projections on SOCs.

\subsection{Tree structure exploitation}
\begin{figure}
    \centering
    \ifarxiv
        \includegraphics{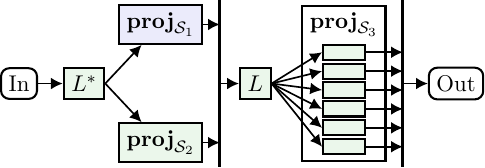}
    \else
        \input{\pathtofigs/forkJoinMain.tex}
    \fi
    \caption{Fork-join chain of streams of main parallelized processes in a CP iteration.
    }
    \label{fig:fork-join-main}
\end{figure}

In order to extract the most benefit from parallel computing, we must first understand the parallel nature of each process in the solver.
Figure~\ref{fig:fork-join-main} illustrates the fork-join chain of the main parallelized processes in a CP iteration.
The projection on the dynamics ($\S_{1}$) is parallelized per stage.
The projections on sets $\S_{2}$ and $\S_{3}$, and the operators $L$ and $L^{\ast}$, are all parallelized per node.
However, there is serial summations across siblings within the dynamics projection and $L^{\ast}$ operator.
As such, a tree structure with high initial branching and no further branching will benefit the most from parallelization.

\subsection{Anderson's acceleration}\label{sec:aa-parallel}
In line~\ref{line:aa} of Algorithm~\ref{alg:spock-gpu}, we follow the same method as before~\cite[Sec. 6.2]{bodard2023spock}, with memory length $m=3$, to compute the directions required. Here, however, we use the CUDA libraries cuSOLVER and cuBLAS to compute the least squares solution. In particular, we compute the QR factorization of a tall matrix with \texttt{cusolverDn<t>geqrf} that stores the resulting orthogonal matrix $Q$ as a sequence of Householder vectors. Then, we use \texttt{cusolverDn<t>ormqr} and \texttt{cublas<t>trsm} to solve the triangular linear system. With this approach, we compared longer memory lengths of up to $m=10$, and found that $m=3$ performed best in all test cases.

This approach can be improved by initially computing the QR factorization in this way, and then updating the Q and R matrices using Givens rotations~\cite{hsieh1993givens}.

\subsection{GPUtils}
Much of \SPOCK{} is built on the lightweight \texttt{GPUtils}~\cite{gputils}, an open-source header-only library. \texttt{GPUtils} makes managing and manipulating device data easy and memory-safe. In particular, this library brings easy \texttt{C++}-type linear algebra operators to \texttt{CUDA C++}.

\section{Results}
As the CP method shares many similarities with ADMM~\cite[Sec. 4.3]{chambolle_method}, our first experiment compared the popular ADMM algorithm against the CP method with the same parallelized code, however, this resulted in more iterations, and more computation per iteration, to obtain the same tolerance.

\begin{figure}
    \centering
    \ifarxiv
        \includegraphics{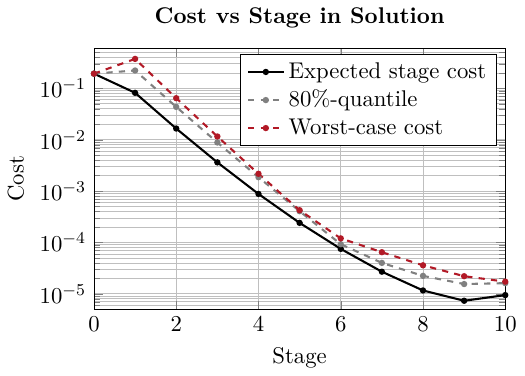}
    \else
        \input{\pathtofigs/costError.tex}
    \fi
    \caption{Cost vs Stage plot for an instance of the problem in~\cite[Sec. 7]{bodard2023spock} with a horizon of 10.
    }
    \label{fig:stage-wise-costs}
\end{figure}

The next experiment gave an insight into the solution returned by the solver. Figure~\ref{fig:stage-wise-costs} shows the stage-wise distribution of costs for an instance of the problem in~\cite[Sec. 7]{bodard2023spock}, with $n_x,n_u,n_w,N=10$ and $n_b=1$. As we expect, the expected cost goes to zero as the stage progresses.

\subsection{Speedup compared to serial \SPOCK{}}
\begin{figure}
    \centering
    \ifarxiv
        \includegraphics{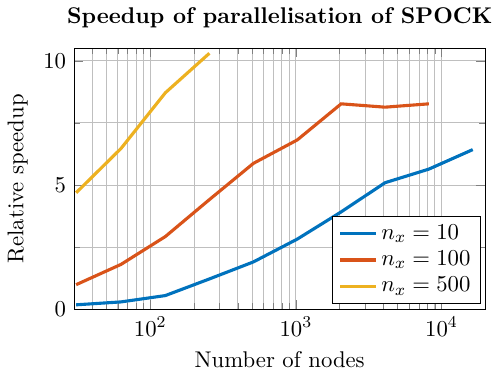}
    \else
        \input{\pathtofigs/speedup_cp.tex}
    \fi
    \caption{Speedup factor of GPU-accelerated \SPOCK{} compared to serial \SPOCK{} for the problem in~\cite[Sec. 7]{bodard2023spock}, where $n_x$ is the number of states and also the number of inputs, i.e., $n_u = n_x$.}
    \label{fig:speedup-cpu-gpu}
\end{figure}

To evaluate the performance improvement of GPU-accelerated \SPOCK{}, we compare its execution time of the problem in~\cite[Sec. 7]{bodard2023spock} against the serial implementation for different problem sizes. Figure \ref{fig:speedup-cpu-gpu} illustrates the speedup factor, defined as the ratio of the execution time of the serial \SPOCK{} to that of the GPU-accelerated \SPOCK.

The GPU-accelerated \SPOCK{} achieves significant speedup across larger problem sizes, both in terms of number of nodes and number of states and inputs. This demonstrates the effectiveness of leveraging parallel computing hardware to handle the computationally intensive tasks in the \SPOCK{} algorithm.

\section{Applications}

\subsection{Case study 1 (Benchmarks)}
To evaluate the solve time of \SPOCK, we compare it against the popular solvers \MOSEK, \GUROBI, and \IPOPT{} for a variety of random problems.

\subsubsection{Problem statement}
We generated a random collection of (large to very large)
risk-averse optimal control problems
(cf. Equation  \eqref{eq:disturbed-dt-affine-system})
along the lines of~\cite[Sec. A.3]{stellato2020osqp}.
We sample prediction horizons $N{}\sim{} \mc{U}(5, 15)$,
and jointly sample $n_b$, $n_w$, and $n_u$
from the uniform distribution on the set
\begin{equation*}
    \left\{
    \begin{bmatrix}
        b
        \\
        w
        \\
        u
    \end{bmatrix}
    \in\N^{3}
    \left|
    \begin{array}{l}
        1 \leq b \leq 3,
        \\
        2 \leq w \leq 10,
        \\
        10 \leq u \leq 300,
        \\
        10^{3} \leq n_v \leq 10^{5}
    \end{array}
    \right.
    \right\},
\end{equation*}
where $n_v$ is the number of variables, i.e.,
${n_v =n_x\cdot|\nodes(0, N)|+n_u\cdot|\nodes(0, N-1)|}$.
We fix $n_x {}={} 2 n_u$.
We impose box state and input constraints,
$-\bar{x} \leq x_{t+1} \leq \bar{x}$,
and
$-\bar{u} \leq u_t \leq \bar{u}$,
for all $t\in\N_{[0, N-1]}$,
with $\bar{x} \sim \mc{U}([1, 2]^{n_x})$
and $\bar{u} \sim \mc{U}([0, 0.1]^{n_u})$.
We use $\AVaR_{\gamma}$ 
as the risk measure with parameter $\gamma\sim\mc{U}([0, 1])$
with a random base probability vector 
$\pi {}\sim{} \mc{U}(\Delta_{n_w})$, where $\Delta_{n_w}$ 
is the probability simplex in $\R^{n_w}$.
Regarding the system dynamics, we take a random 
matrix $B{}\in{}\R^{n_x \times n_u}$ with 
elements $B_{ij}\sim\mc{N}(0, 1)$ and 
for each event $w\in\N_{[1, n_w]}$ we take 
random system matrices 
\begin{subequations}%
\begin{align}
    \label{eq:random-system-dynamics-A}%
    A(w)
    {}={} 
    I_{n_x} + A^w, A^w_{ij}\sim\mc{N}(0, 0.01),
    \\
    \label{eq:random-system-dynamics-B}
    B(w)
    {}={} 
    B + B^w, B^w_{ij}\sim\mc{N}(0, 0.01).
\end{align}%
For the stage costs we define 
\begin{align}
    Q_0 = \diag{q_0}, q_0\sim\mc{U}([0, 0.1]^{n_x}),
    \\
    R_0 = \diag{r_0}, r_0\sim\mc{U}([0, 100]^{n_u}),
\end{align}
and for each event $w\in\N_{[1, n_w]}$
\begin{align}
    Q^w_{ij}{}\sim{}&\mc{N}(0, 0.01),
    \\
    Q(w) {}={}& (Q_0 + Q^w)(Q_0 + Q^w)\Tr,
    \\
    R^w_{ij}{}\sim{}&\mc{N}(0, 0.01),
    \\
    R(w) {}={}& (R_0 + R^w)(R_0 + R^w)\Tr, 
\end{align}
\end{subequations}
and $Q_N = Q_0$. 
Lastly, we sample $x\in\R^{n_x}$ with 
$x_i \sim\mc{U}([-0.5\bar{x}_i, 0.5\bar{x}_i])$.

\subsubsection{Performance profiles}
We use performance profiles \cite{dolan2002benchmarking} to compare the performance of the \SPOCK{} algorithm with state-of-the-art solvers. We define the performance measure to be the solve time. We limit the solve time to \SI{5}{\minute}, and we consider exceeding this limit a failure. Failure is given the ratio $r_{\rm M} = \infty$.

\begin{figure}
    \centering
    \ifarxiv
        \includegraphics{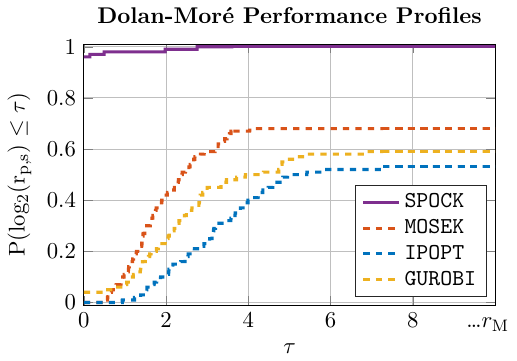}
    \else
        \input{\pathtofigs/dolanMore_random.tex}
    \fi
    \caption{Dolan-Mor\'e plot of popular solvers \MOSEK, \GUROBI, and \IPOPT{} compared to GPU-accelerated \SPOCK{} for 100 random problems.}
    \label{fig:dolan-more-random}
\end{figure}

We present the Dolan-Mor\'{e} plot of the performance profiles in Figure~\ref{fig:dolan-more-random}. The plot shows the cumulative distribution of the performance ratios for the solvers \MOSEK, \GUROBI, \IPOPT, and GPU-accelerated \SPOCK{} for a set of $100$ random problems. The plot shows the percentage of problems solved within a factor $\tau$ of the best solver. The best solver is the one with the shortest solve time for each problem. The plot shows that \SPOCK{} outperforms the other solvers in terms of speed (highest probability at $\tau=0$) and robustness (highest probability at $\tau=r_{\rm M}$).

\subsubsection{Memory management}
\begin{figure}
    \centering
    \ifarxiv
        \includegraphics{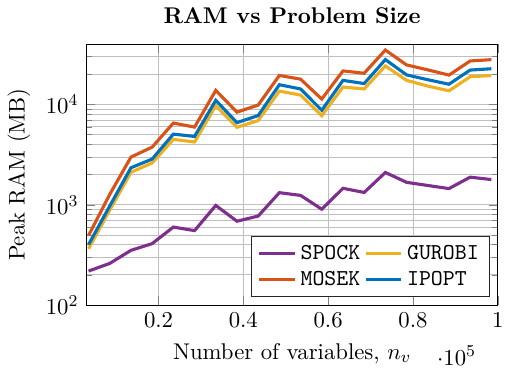}
    \else
        \input{\pathtofigs/ramUsage_random.tex}
    \fi
    \caption{Peak RAM vs problem size of solvers \MOSEK, \GUROBI, \IPOPT, and GPU-accelerated \SPOCK{} for 100 random problems.}
    \label{fig:peak-ram-random}
\end{figure}

The peak RAM usage of each solver was recorded while solving the 100 random problems. CPU and GPU memory usage were combined for GPU-accelerated SPOCK{}, and the others only use CPU memory. The results are presented in Figure~\ref{fig:peak-ram-random}, illustrating how \SPOCK{} uses significantly less memory than the other solvers.
For example, for a problem size of $n_v=10^5$, \SPOCK{} uses around \SI{1}{\giga\byte} of memory, while the other solvers 
use around 25--30 \SI{}{\giga\byte} of memory.

\subsection{Case study 2 (Parallelization)}
Here we show the effect of different tree branching on the parallelization, and, therefore, the solve times, of the algorithm.
The problem statement is adapted from~\cite{patrinos2011networked}.
Consider a networked control system (NCS) 
consisting of a continuous-time LTI plant and a discrete-time controller that
are connected through a communication network with induced 
sensor-to-controller (SC) delay, $\sigma$~\cite{cloosterman2008control}.
The full
state of the system, $x$, is sampled by a sensor
with a constant sampling interval $T > 0$. 
The discrete-time controller is event-driven and able to monitor the SC delay, via timestamping. 
The discrete-time control signal $u_k$ is transformed to a continuous-time control input $u(t)$ by a zero-order hold.

\subsubsection{Sensor-to-controller delay}
The SC delay is modeled as a discrete-time stochastic process $\sigma$ with a.s. bounded paths in $[0, T]$ and the martingale property 
$\E[\sigma_{t+1}|\sigma_{t}] = \sigma_{t}$.

First, we consider a stochastic process, $(\phi_{t})_t$, that takes values in $[0, 1]$. We assume $\phi_{t+1}\mid \phi_t \sim {\rm Beta}(\zeta_t, \beta_t)$ and we choose $\zeta_t = \zeta_t(\phi_t)$ and $\beta_t = \beta_t(\phi_t)$ so that $\E[\phi_{t+1}\mid \phi_t] = \phi_t$ is satisfied; then, we define $\sigma_{t} = \phi_t T$.

We impose
$
\E[\phi_{t+1}\mid \phi_t]
= 
\frac{\zeta_t}{\zeta_t+ \beta_t} = \phi_t,
$
or, equivalently,
\begin{align}
\frac{\zeta_t}{\beta_t} = \frac{\phi_t}{1-\phi_t}.
\end{align}
To achieve this we assume $\zeta_t = \theta \phi_t$ and $\beta_t = \theta (1-\phi_t)$, so that the above is automatically satisfied, and $\theta > 0$ is a free parameter. It can then be seen that 
${\rm Var}[\phi_{t+1}{}\mid{}\phi_t]=\frac{\phi_t(1-\phi_t)}{\theta+1}.$
This means that a large $\theta$ leads to a small conditional variance.
The unconditional variance is 
${\rm Var}[\phi_t] = \phi_0(1-\phi_0)(1 - \tilde{\theta}^t),$
where $\tilde{\theta} = \frac{\theta}{\theta+1}$, so ${\rm Var}[\phi_t]\to \phi_0(1-\phi_0)$ as $t\to\infty$.

\begin{figure}
    \centering
    \ifarxiv
        \includegraphics{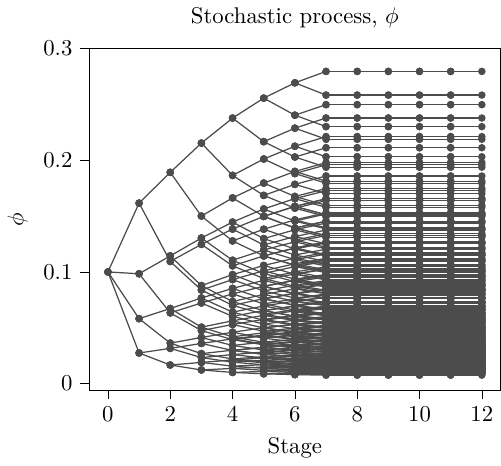}
    \else
        \input{\pathtofigs/delay_tree.tex}
    \fi
    \caption{Example scenario tree of stochastic process, $\phi$, where $\theta$=1 and $\phi^{0}$=0.1.}
    \label{fig:time-delay-scenario-tree}
\end{figure}

Now, we construct a scenario tree by sampling $\phi$, where the children of a node are equiprobable.
In particular, having determined the value of $\phi^i$ of a node $i\in\nodes(0, N-1)$, we sample $\phi^{i_+}\sim{\rm Beta}(\theta \phi^i, \theta(1-\phi^i)), \,\forall i_+\in\ch(i)$.
That is, we split $[0,1]$ into $|\ch(i)|$ many intervals of equal length, and take their centers $\iota^{i_+}$. Then, $\phi^{i_+}$ is the value of the quantile of the above Beta at $\iota^{i_+}$.
For this case study we set $\phi^{0}=0.1$ and $\theta=1$. An example of the scenario tree is shown in Figure~\ref{fig:time-delay-scenario-tree}.

\subsubsection{System dynamics and constraints}
The NCS model is
\begin{subequations}\label{eq:networked-system-dynamics}
\begin{align}
    \dot{x}(t) 
    {}={}& 
    A_c x(t) + B_c u(t), 
    \\
    u(t) {}={}& u_k, 
    \quad 
    t \in 
    \left[\left.
        kT + 
        \sigma_k,
        \ 
        (k+1)T + 
        \sigma_{k+1}
    \right)\right.,
\end{align}
\end{subequations}
where 
\begin{subequations}
\begin{align}
    A_c 
    {}={}&
    I_{n_x} + \diag{a}, \quad
    a_{i}\overset{\rm iid}{\sim}\mc{N}(0, 0.01),
    \\
    B_c 
    {}={}&
    0.1 \cdot 1_{n_x \times n_u}.
\end{align}
\end{subequations}
Using the technique described in~\cite{patrinos2011networked},
\eqref{eq:networked-system-dynamics} is transformed into a discrete-time system of the form in \eqref{eq:finite-horizon-evolution-2-gpu}
by letting the augmented system state be
$
\tilde{x}_k = 
\begin{bmatrix}
    x_{k}\Tr & u_{k-1}\Tr
\end{bmatrix}\Tr 
\in \R^{n_x + n_u}
$
and $x_k = x(kT)$.
That is,
\begin{equation}
    \tilde{x}^{i} = A^{i} \tilde{x}^{i_-} + B^{i} u^{i_-}
\end{equation}
where $A^{i}$ and $B^{i}$ are functions of $A_c$, $B_c$, and $\sigma^{i}$, following from \cite[Eq. (5)]{patrinos2011networked}.
The system is subject to discrete-time state, 
$\|x^{i}\|_{\infty} \leq 3, \forall i\in\nodes(0,N)$,
and input constraints
$\|u^{i}\|_{\infty} \leq 0.9, \forall i\in\nodes(0,N-1)$.

\subsubsection{Cost}
We use a quadratic stage cost function
\begin{equation}
    \ell(\tilde{x}^i, u^i)
    {}={}
    \tilde{x}\iTr Q \tilde{x}^i + u\iTr R u^i,
\end{equation}
where 
$Q = \diag{10^{-6} \cdot 1_{n_x+n_u}}$ 
and 
$R = \diag{0.1 \cdot 1_{n_u}}$.
We use a quadratic terminal cost function
\begin{equation}
    \ell_N(\tilde{x}^j)
    {}={}
    \tilde{x}\jTr Q_N \tilde{x}^j,
\end{equation}
where 
$Q_N = \diag{
    \begin{bmatrix}
        1_{n_x}\Tr & 
        10^{-6} \cdot 1_{n_u}\Tr
    \end{bmatrix}}
$.

\subsubsection{Performance}
Here we investigate the effect of different tree branching on the parallelization of the algorithm, and so on the solve time. We also demonstrate the speedup of \SPOCK{} over current state-of-the-art solvers.
We set the horizon to $N=48$, $n_x=100$, $n_u=50$, and the sampling time is $T=5$.
The risk measure is $\AVaR_{0.95}$.
The initial state and previous input are set element-wise to $2$ and $0$, respectively.
\begin{figure}
    \centering
    \ifarxiv
        \includegraphics{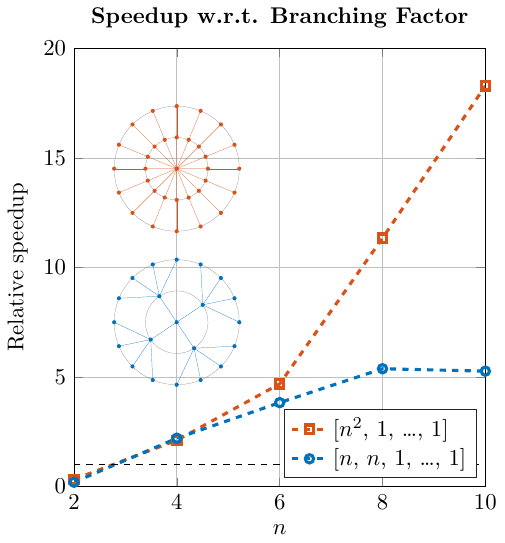}
    \else
        \input{\pathtofigs/delay_speedup.tex}
    \fi
    \caption{Speedup against problem size of GPU-accelerated \SPOCK{} compared to the next fastest solver between \MOSEK, \GUROBI, \IPOPT. We consider two distinct stage branching factors with the same number of scenarios, illustrated for the first two stages at $n=4$.}
    \label{fig:speedup-power}
\end{figure}
Figure~\ref{fig:speedup-power} demonstrates that a tree with a high initial branching (stage branching of $[n^2, 1, \ldots, 1]$) is more amenable to parallelization than one with further branching (stage branching of $[n, n, 1, \ldots, 1]$). Note that the algorithm step size $\alpha$ is the same for both cases due to preconditioning.
Figure~\ref{fig:speedup-power} further demonstrates the capability of \SPOCK{} over \MOSEK, \GUROBI, and \IPOPT{} with increasing speedup for wider problems (i.e., more scenarios).

\section{Conclusions}
This work demonstrates the feasibility and advantages of GPU acceleration for solving RAOCPs. GPU-accelerated \SPOCK{} significantly outperforms its serial counterpart and current state-of-the-art interior-point solvers, making it a strong candidate for real-time control applications with uncertain systems.

\section*{Acknowledgements}
The authors gratefully acknowledge the invaluable guidance and support of Christian A. Hans, with Automation and Sensorics
in Networked Systems Group, University of Kassel, Germany.

\printbibliography[title={References},heading=refs]

\end{document}